\documentclass[11pt]{article}
\usepackage{latexsym}
\usepackage{amsmath}
\usepackage[dvips]{graphicx}
\usepackage{amsfonts}
\usepackage{amstext}
\usepackage{amssymb}
\usepackage{amsxtra}

\newtheorem{lemma}{Lemma}[section]
\newtheorem{theorem}{Theorem}[section]
\newtheorem{remark}{Remark}[section]

\def\IR{\rm I \kern-0.20em R}

\newcommand{\EE}{{\mathord{I\kern -.33em E}}}
\newcommand{\PP}{{\mathord{I\kern -.33em P}}}
\newcommand{\vol}{volatility }

\def \BE{\begin{equation}}
\def \en{\end{equation}}

\title{Estimating the Fractal Dimension of the S\&P 500 Index using Wavelet Analysis}
\author{Erhan Bayraktar \thanks{Department of Mathematics, University of Michigan, 2074 East Hall, Ann Arbor, MI 48109,
{\em erhan@umich.edu}} \and H. Vincent Poor
\thanks{Department of Electrical Engineering, Princeton University
, E-Quad, Princeton, NJ 08544, {\em poor@princeton.edu}}\and K.
Ronnie Sircar
\thanks{Department of Operations Research \& Financial
Engineering, Princeton University , E-Quad, Princeton, NJ 08544,
{\em sircar@princeton.edu}}} 
\begin{document}
\maketitle

\begin{abstract}
S\&P 500 index data sampled at one-minute intervals over the
course of 11.5 years (January 1989- May 2000) is analyzed, and in
particular the Hurst parameter over segments of stationarity (the
time period over which the Hurst parameter is almost constant) is
estimated. An asymptotically unbiased and efficient estimator
using the log-scale spectrum is employed. The estimator is
asymptotically Gaussian and the variance of the estimate that is
obtained from a data segment of $N$ points is of order
$\frac{1}{N}$. Wavelet analysis is tailor made for the high
frequency data set, since it has low computational complexity due
to the pyramidal algorithm for computing the detail coefficients.
This estimator is robust to additive non-stationarities, and here
it is shown to exhibit some degree of robustness to multiplicative
non-stationarities, such as seasonalities and volatility
persistence, as well. This analysis shows that the market became
more efficient in the period 1997-2000.


\end{abstract}
 \footnotetext[2]{$\textbf{Key Words:}$
High-frequency data, S\&P 500 index, long range dependence, heavy
tailed marginals, fractional Brownian motion, wavelet analysis,
log scale spectrum}

\newpage
\section{Introduction\label{introduction}}
Stochastic models based primarily on continuous or discrete time
random walks have been the foundation of financial engineering
since they were introduced in the economics literature in the
1960s. Such models exploded in popularity because of the
successful option pricing theory built around them by Black and
Scholes
 \cite{kn:blackscholes} and Cox {\em et al.}  \cite{kn:crr}, as
well as the simplicity of the solution of associated optimal
investment problems given by Merton  \cite{kn:merton69}.

Typically, models used in finance are diffusions built on standard
Brownian motion and they are associated with partial differential
equations describing corresponding optimal investment or pricing
strategies. At the same time, the failure of models based on
independent increments to describe certain financial data has been
observed
since Greene and Fielitz \cite{kn:fielitz} and Mandelbrot
\cite{kn:mandel}, and \cite{kn:mandellong}. Using R/S analysis,
Greene and Fielitz studied 200 daily stock returns of securities
listed on the New York Stock Exchange and they found significant
long range dependence. Contrary to their finding, Lo \cite{kn:lo},
using a modified R/S analysis designed to compensate for the
presence of short-range dependence, finds no evidence of
long-range dependence (LRD). However, Teverovsky et al.
\cite{kn:taqqu2} and Willinger et al. \cite{kn:taqqu1} identified
a number of problems associated with Lo's method. In particular,
they showed that Lo's method has a strong preference for accepting
the null hypothesis of no long range dependence. This happens even
with long-range dependent synthetic data.
To account for the long-range dependence observed in financial
data Cutland et al. \cite{kn:cutland} proposed to replace Brownian
motion with {\em fractional Brownian motion} (fBm) as the building
block of stochastic models for asset prices. An account of the
historical development of these ideas can be traced from Cutland
et al \cite{kn:cutland}, Mandelbrot \cite{kn:mandelbrot} and
Shiryaev \cite{kn:shiryaev}. The S\&P 500 index was analyzed in
\cite{kn:peters1} and \cite{kn:peters} by Peters using R/S
analysis, and he concluded that the raw return series exhibits
long-range dependence. See also \cite{kleinow} for analysis of LRD
in German stock indices.

Here we present a study of a high-frequency financial data set
exhibiting long-range dependence, and develop wavelet based
techniques for its analysis. In particular we examine the S\&P 500
over 11.5 years, taken at one-minute intervals. The wavelet tool
we consider, namely the log-scale spectrum method, is
asymptotically unbiased and efficient with a vanishing precision
error for estimating the Hurst parameter (a measure of long-range
dependence, explained in ~\eqref{eq:covfbm} below). (See Theorem
\ref{mainthm2}.) Since we are dealing with high frequency data, we
need fast algorithms for the processing of the data. Wavelet
analysis is tailor-made for this purpose due to the pyramidal
algorithm, which calculates the wavelet coefficients using octave
filter banks. In essence, we look at a linear transform of the
logarithm of the wavelet variance ({\em i.e.} the variance of the
detail coefficients, defined in~\eqref{eq:detail}) to estimate the
Hurst parameter.
Moreover, the log-scale spectrum methodology is insensitive to
additive non-stationarities, and, as we shall see, it also
exhibits robustness to multiplicative non-stationarities of a very
general type including seasonalities and volatility persistence
(Section \ref{presistence}).

Although the Hurst parameter of S\&P 500 data considered here is
significantly above the efficient market value of $H=\frac{1}{2}$,
it began to approach that level around 1997. This behavior of the
market might be related to the increase in Internet trading, which
has the three-fold effect of increasing the number of small
traders, increasing the frequency of trading activity, and
improving traders' access to price information. An analytical
model of this observation is proposed in \cite{kn:bhs}.


\subsection{Fractional Brownian Motion} A natural extension of
the conventional stochastic models for security prices to
incorporate long-range dependence is to model the price series
with geometric fractional Brownian motion:
\begin{equation}\label{eq:model}
P_{t}=P_{0}\exp\left(\mu_{t}+\int_{0}^{t}\sigma_{s}
dB^{H}_{s}\right),
\end{equation}
where $P_0$ is today's observed price, $\mu$ is a  growth rate
parameter, $\sigma$ is the stochastic \vol process, and $B^{H}$ is
a fractional Brownian motion, an almost surely (a.s.) continuous
and centered Gaussian process with stationary increments.
autocorrelation of $B^H$
\begin{equation}
\EE\left\{B^{H}_{t}B^{H}_{s}\right\}=\frac{1}{2}\left(|t|^{2H}+|s|^{2H}-|t-s|^{2H}\right),~\label{eq:covfbm}
\end{equation}
where $H \in (0,1]$ is the so-called Hurst parameter. (Note that
$H=\frac{1}{2}$ gives standard Brownian motion.) From this
definition, it is easy to see that fBm is self-similar, {\em i.e.}
$B^{H}(at)=a^{H}B(t)$, where the equality is in the sense of
finite dimensional distributions. This model for stock market
prices is a generalization of the model proposed in
\cite{kn:cutland} to allow for non-Gaussian returns distribution
into the model. Heavy tailed marginals for stock price returns
have been observed in many empirical studies since the early
1960's by Fama \cite{kn:fama} and Mandelbrot
\cite{kn:mandelbrot2}.

Fractional Brownian motion models are able to capture long range
dependence in a parsimonious way. Consider for example the
fractional Gaussian noise $Z(k):=B^{H}(k)-B^{H}(k-1)$. The
auto-correlation function of $Z$, which is denoted by $r$,
satisfies the asymptotic relation
\begin{equation}
r(k) \sim r(0)H(2H-1)k^{2H-2}, \quad \text{as $k \rightarrow
\infty$}.
\end{equation}
For $H \in (1/2,1]$, $Z$ exhibits long-range dependence, which is
also called the Joseph effect in Mandelbrot's terminology
\cite{kn:mandelbrot}. For $H=1/2$ all correlations at non-zero
lags are zero. For $H \in (0,1/2)$ the correlations are summable,
and in fact they sum up to zero. The latter case is less
interesting for financial applications (\cite{kn:cutland}).



Now, we will make the meaning of (\ref{eq:model}) clear by
defining the integral term. The stochastic integral in
(\ref{eq:model}) is understood as the probabilistic limits of
Stieltjes sums. That is, given stochastic processes $Y$ and $X$,
such that $Y$ is adapted to the filtration generated by $X$, we
say that the integral $\int YdX$ exists if, for every $t<\infty$,
and for each sequence of partitions $ \{\sigma^{n}\}_{n \in
\mathbb{N}}$, $\sigma^{n}=(T^n_{1},T^n_{2},...,T^n_{k_{n}})$, of
the interval $[0,t]$ that satisfies $\lim_{n \rightarrow \infty},
\max_{i} |T^n_{i+1}-T^n_{i}|=0$, the sequence of sums
$\left(\sum_{i}Y_{T^{n}_{i}}(X_{T^n_{i+1}}-X_{T^n_{i}})\right)$
converges in probability. That is, we define
\begin{equation}\label{stieltjes}
\int_{0}^{t}Y_{s}dX_{s}=\mathbb{P}-\lim_{n \rightarrow
\infty}\sum_{i}Y_{T^{n}_{i}}(X_{T^n_{i+1}}-X_{T^n_{i}}).
\end{equation}
By the Bichteler-Dellacherie Theorem \cite{kn:Pro} one can see
that the integrals of adapted processes with respect to fBm may
not converge in probability. However when $H>\frac{1}{2}$ there
are two families of processes that are integrable with respect to
fBm that are sufficiently large for modeling purposes. The first
family consists of continuous semi-martingales adapted to the
filtration of fBm  as demonstrated in \cite{kn:erhan}. The second
family consists of processes with H\"{o}lder exponents greater
than $1-H$. (This integration can be carried out pathwise as
demonstrated in \cite{kn:ruzmaikina}, and \cite{kn:zahle}).



\subsection{Markets with Arbitrage Opportunities}
Much of finance theory relies on the assumption that markets
adjust prices rapidly to exclude any arbitrage opportunities. It
is well known that models based on fBm allow arbitrage
opportunities  (\cite{kn:cheridito} and \cite{kn:rogers}). Even in
the case of stochastic $\sigma$ we have shown that there exist
arbitrage
opportunities in a single stock setting  \cite{kn:erhan}. 
However, strategies that capitalize on the smoothness (relative to
standard Bm) and correlation structure of fBm to make gains with
no risk, involve exploiting the fine-scale properties of the
process' trajectories. Therefore, this kind of model describes a
market where arbitrage opportunities can be realized (by frequent
trading), which seems plausible in real markets. But the ability
of a trader to implement this type of strategy is likely to be
hindered by market frictions, such as transaction costs and the
minimal amount of time between two consecutive transactions.
Indeed Cheridito \cite{kn:cheridito} showed that by introducing a
minimal amount of time $h>0$ between any two consecutive
transactions, arbitrage opportunities can be excluded from a
geometric fractional Brownian motion model (\emph{i.e} when
$\sigma$ is taken to be constant in (\ref{eq:model})).


Elliot and Van der Hoek \cite{elliot}, and Oksendal and Hu
\cite{kn:oksendal2} considered another fractional Black-Scholes
(B-S) model by defining the integrals in (\ref{eq:model}) as Wick
type integrals. This fractional B-S model does not lead to
arbitrage opportunities; however one can argue that it is not a
suitable model for stock price dynamics. The Wick type integral of
a process $Y$ with respect to a process $X$ is defined as
\begin{equation}\label{eq:wick}
\int_{0}^{t}Y_{T^n_{i}}\diamond(X_{T^n_{i+1}}-X_{T^n_{i}})
\end{equation}
where the convergence is in the $L^2$ space of random variables.
(The Wick product is defined using the tensor product structure of
$L^2$; see \cite{kn:hida}.) The Wick type integral of $Y$ with
respect to fBm with Hurst parameter $H$ is equal to the Stieltjes
integral defined above plus a drift term (see \cite{kn:duncan}
Thm. 3.12),
\begin{equation*}
\int_{0}^{t}Y_{s}dB^{H}_{s}=\int_{0}^{t}Y_{s}\diamond
dB^{H}_{s}+\int_{0}^{t}D^\phi_{s}Y_{s}ds,
\end{equation*}
where $\phi(s,t)=H(2H-1)|s-t|^{2H-2}$, and $D^\phi_{s}
Y_t:=(D^\phi Y_t)(s)$ is the Hida derivative of the random
variable $Y_t$. Hence writing an integral equation in terms of
Wick product integrals is equivalent to writing a Stieltjes
differential equation with a different drift term.  The fractional
B-S model with the integrals defined as in (\ref{eq:wick}) does
not lead to arbitrage opportunities.  However, this conclusion is
based on the redefinition of the class of self-financing
strategies. The self-financing strategies in a Stieltjes framework
are no longer self-financing strategies in a Wick framework, so
that all the self-financing arbitrage strategies of the Stieltjes
framework are ruled out by the approach of \cite{elliot} and
\cite{kn:oksendal2}. However in the Wick framework it is hard to
give economic interpretations to trading strategies. For
illustration let us consider a simple hold strategy. Let $u$
denote the number of shares that are held at time $T_{1}$ by an
economic agent, and let us see the value change of the portfolio
over the time interval $[T_{1},T_{2}]$ if the agent chooses to
hold its shares of the risky asset in a Wick type framework. If
$P_{t}$ denotes the price of the risky asset at time $t$, the
increment of the value of the portfolio over the interval
$[T_{1},T_{2})$ is
\begin{equation}\label{wick}
u \diamond (P_{T_{2}}-P_{T_{1}}).
\end{equation}
It is hard to attach a clear economic meaning to this quantity
since the Wick product is not a path-wise product but rather is
defined using the tensor product structure of the $L^2$ space of
random variables. On the other hand, (\ref{stieltjes}) involves
the actual realization of the increment,
\begin{equation}\label{stielt}
u(\omega) (P_{T_{2}}(\omega)-P_{T_{1}}(\omega)),
\end{equation}
which has a clear economic interpretation. (Here, $\omega$ denotes
the point in the sample space corresponding to the given
realization of the price process.) So, among the two candidates
for the value of the increment of a simple hold strategy,
(\ref{stielt}) has a more direct economic meaning. Hence the no
arbitrage conclusion of \cite{elliot} and \cite{kn:oksendal2}
cannot be interpreted within the usual meaning of this term, and
thus we prefer to apply the definition (\ref{stieltjes}) for the
stochastic integrals involved. (Also see \cite{bjork} and
\cite{valkeila} which also argue that Wick type integrals are not
suitable for defining trading strategies.)

Models with fBm differentials in the stochastic differential
equations describing the stock price can be built however, by
considering the Nash-equilibrium which arises from a game in which
the players are instituional investors manipulating the
coefficients of a stochastic differential equation with fBm
differentials in order to maximize their utilities. The
Nash-equilibrium for such stochastic differential games is
considered by Bayraktar and Poor in \cite{kn:bp}. The fBm
differentials in the controlled stochastic differential game can
be interpreted as the trading noise arising from the activities of
small investors who exhibit inertia (see \cite{kn:bhs}).

\subsection{Non-stationarities Expected from a Financial Time
Series}

\subsubsection{Time-Variation of $H$} In this paper,
we are interested in the estimation of the Hurst parameter ($H$)
from historical stock index data. In addition, we will study the
variability of this parameter over time. Common experience with
financial data suggests that it exhibits too much complexity to be
described by a model as simple as (\ref{eq:model}), which says
that the log price process
\begin{equation}\label{logstock}
Y_t := \log (P_t/P_0)=\mu_{t}+\int_{0}^{t}\sigma_{s} dB^{H}_{s},
\qquad t \in [0,T],
\end{equation}
 is a stochastic integral with respect to fBm
with drift. In particular, if we could remove the
non-stationarities due to the drift and stochastic volatility,
then $Y_{t}$ would be a process with stationary increments.
However, stationarity is not usually a property of return series
of a financial index, which are often extremely turbulent.
Therefore, we would like to identify segments of time over which
the return series is close to stationary. In other words, one of
our aims is to study the variation of $H$ over time as a gauge of
the epochs when the returns process behaves like a stationary
process. We do not assume any particular form of temporal behavior
for this parameter; its variation is to be found from our data
analysis. We partition the data into smaller segments, find the
corresponding parameters for each of the segments, and use
filtering to remove the extrinsic variation in the parameters due
to finiteness of the segment.  We vary the segmentation size and
repeat the procedure described. Then, comparing the fluctuations
of $H$ among different segmentation levels, we are able to come up
with the segments of stationarity. The comparison among the
different levels of segmentation is possible since extra noise
introduced by altering the segmentation level is filtered out.
In Section 5, we show how we come to the conclusion that the
segments of stationarity for the S\&P 500 are on the order of
$2^{14}$ points, or approximately $8$ weeks. 

\subsubsection{Drift and Stochastic Volatility}

We can also allow the average growth rate $\mu$ in
(\ref{eq:model}) to have time variation. Our analysis is
insensitive to polynomial trends of certain order. Our method,
based on the analysis of the log-scale spectrum, is also
insensitive to additive periodic components. The effects of
periodicity on the scale spectrum, and a technique for alleviating
the polluting effects of additive periodicity by increasing the
number of vanishing moments of a mother wavelet is analyzed by
Abry \emph{et al.} \cite{kn:veitch1} on fBm with $H=0.8$ and an
additive sinusoidal trend.

Here we are interested in S\&P 500 data, for which the returns
have a multiplicative periodic component and stochastic volatility
in addition to their intrinsic random variation. The existence of
seasonalities is observed in various financial time series: see
\cite{kn:barucci} for a single stock return series,
\cite{kn:fouque} for S\&P 500 index data, and \cite{kn:andersen2}
and \cite{kn:andersen3} for FEX data. Heavy tailed marginals for
stock price returns have been observed in many empirical studies
since the early 1960's (e.g.,
\cite{kn:fielitz},\cite{kn:mandelbrot2}). So we expect to have
stochastic volatility\footnote{Observe, however, that equation
(\ref{eq:model}) should not be read in the same way as financial
models driven by standard Brownian motion, since the stochastic
integral $\int\sigma\,dB^H$ is not a martingale. Moreover,
$\sigma^{2}$ is not the quadratic variation process of this
integral, and should not be viewed as the volatility process
literally.} in (\ref{eq:model}) as well. Therefore we must take
into account these non-stationarities in the data while developing
an estimation procedure. One way of dealing with seasonalities is
given in \cite{kn:andersen} and \cite{kn:rambo2}.  Here we show
that the scale spectrum method is quite insensitive to
multiplicative nonstationarities as well as to volatility
persistence.


The remainder of this paper is organized as follows. In Section 2,
we introduce our estimation technique and discuss its statistical
and robustness properties. In Section 3 we apply our technique to
S\&P 500 index data and discuss our observations.

\section{The Log-Scale Spectrum Methodology}

In the Appendix we provide a brief introduction to wavelets
(following the treatment by Mallat  \cite{kn:mallat}), and the
pyramidal algorithm and its initialization. Henceforth we will
assume the notation introduced in Section \ref{wavelet}.

For $(j,k) \in \mathbb{Z}^{2}$ let $d_j(k)$ denote the detail
coefficient at scale $j$ and shift $k$ of a random process $Z$:
\begin{equation}\label{eq:detail}
d_{j}(k)=2^{-j/2}\int_{-\infty}^{+\infty}\psi(2^{-j}t-k)Z(t)\,dt,
\end{equation}
where $\psi$ is any function satisfying the vanishing moments
condition (Appendix A.1.1, (\ref{vmoments})) for some $p \in
\mathbb{Z}_{+}$.

The empirical variance as a function of the scale parameter $j \in
\mathbb{Z}$ is called the scale spectrum and is given by
\begin{equation}
S_{j}=\frac{1}{N/2^{j}}\sum_{k=1}^{N/2^{j}}[d_{j}(k)]^{2}, \qquad
j \leq \log_{2}(N), \label{eq:scale}
\end{equation}
where $N$ is the number of initial approximation coefficients,
i.e. for $j=0$.

\subsection{Scale Spectra and fBm}\label{detfbm}

The detail coefficients of fBm satisfy the following,
\begin{equation}
\EE\{[d_{j}(k)]^{2}\}=K(H)2^{(2H+1)j}, ~\label{eq:exponent}
\end{equation}
where $K(H)$ is given by
\begin{equation}
K(H)=\frac{1-2^{-2H}}{(2H+1)(2H+2)},~\label{eq:kh}
\end{equation}
as given in \cite{kn:abry}, for example. The behavior in
~\eqref{eq:exponent} suggests that the empirical variance of the
sequence $\big(d_{j}(k)\big)_{k \leq N/2^{j}}$ can be used to
estimate the Hurst parameter $H$. The empirical variance of fBm
satisfies
\begin{equation*}
\EE\{S_{j}\}=\EE\{[d_{j}(k)]^{2}\}=K(H)\sigma^{2}2^{(2H+1)j}.
\end{equation*}
We immediately see that the slope of the $\log$ scale spectrum
$\log(S_j)$ yields a simple estimator of the Hurst parameter, but
as we shall see one can do better than this simple estimator.

\subsection{Synthetically Generated fBm and Corresponding Log-scale Spectra}
In this section, we illustrate the behavior of the log-spectrum on
synthetic data. We use the method of Abry and Sellan
\cite{kn:abry} to generate a realization of $2^{19}$ points of fBm
with $H=0.6$. This method uses wavelets for the synthesis of fBm
and requires the specification of the number of scales, which we
choose to be $20$. This simulation method is extremely fast, which
is important for our purposes since we need on the order of 1
million data points to carry out our synthetic analysis.
\begin{figure}[ht]
\begin{center}
\includegraphics[width = 0.8\textwidth,height=6cm]{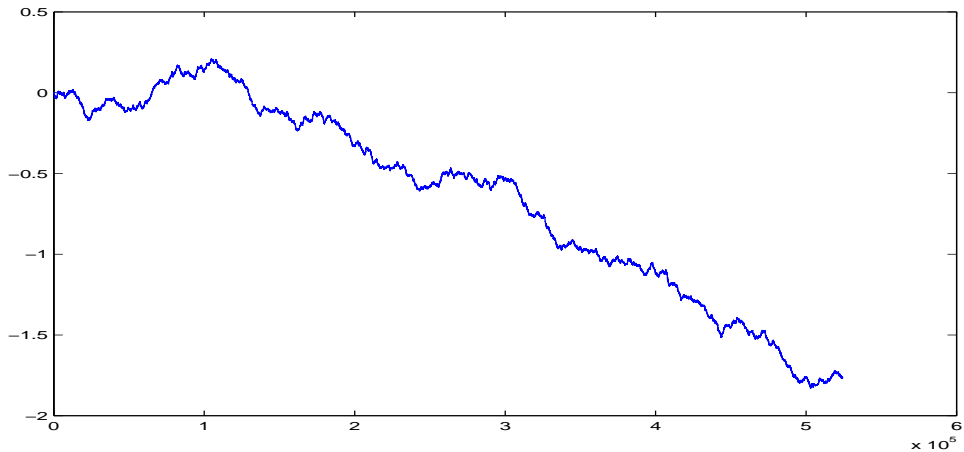}
\caption{A realization of fBm ($H=0.6$) of $2^{19}$ points created
using 20 scales} \label{fig:fbm}
\end{center}
\end{figure}

\begin{figure}
\begin{center}
\includegraphics[width = 0.8\textwidth,height=6cm]{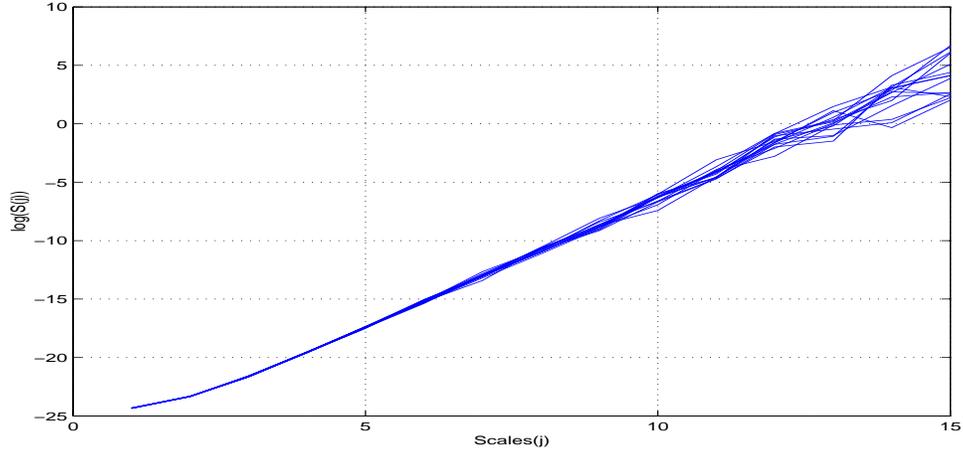}
\caption{Log-scale spectra for segments of length $2^{15}$ for the
realization of fBm in Fig.~\ref{fig:fbm} } \label{fig:fbmscale}
\end{center}
\end{figure}
\begin{figure}
\begin{center}\label{fbmhurts}
\includegraphics[width = 0.8\textwidth,height=6cm]{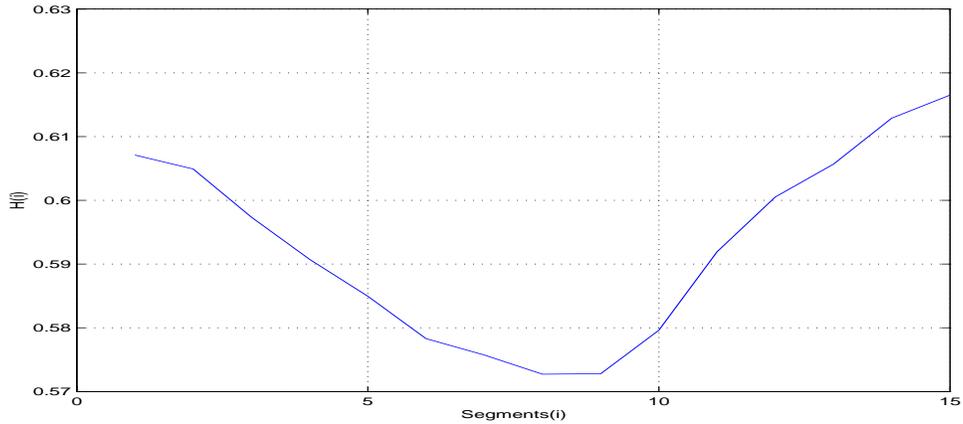}
\caption{Hurst estimates for the realization of fBm given in
Fig.~\ref{fig:fbm} when the segment lengths are $2^{15}$}
\label{fig:fbmh}
\end{center}
\end{figure}

We use segments of length $2^{15}$ (\emph{i.e.}, $2^{15}$
minutes), and the estimates of the Hurst parameter over each
segment for the case of fBm are shown in Fig.~\ref{fig:fbmh}. The
associated log-scale spectra are shown in Fig.~\ref{fig:fbmscale}.
(The mean of the estimates of $H$ over the segments is 0.5928, and
the standard deviation(std) is 0.0149.)

In the next section, we will analyze the asymptotic properties of
the logarithm of the scale spectrum when $\sigma$ in
(\ref{logstock}) is taken to be constant, and we will develop an
asymptotically efficient estimator using these results. Then in
the following section, using the path properties of the integrals
with respect to fBm, we will discuss the robustness of this
estimator to stochastic volatility and seasonalities.

\subsection{Asymptotic Distribution of the Logarithm of the Scale
Spectrum\label{section4}}

We generalize the method developed by Papanicolaou and Solna
\cite{kn:solna} to our case when the process to be analyzed is
given by (\ref{logstock}) with $\sigma$ constant.  The treatment
of \cite{kn:solna} was concerned with Kolmogorov turbulence for
which $H$ is around $\frac{1}{3}$, and therefore the use of a Haar
wavelet in (\ref{eq:detail}) suffices. In the S\&P 500 data, $H$
is expected to be greater than or equal to $\frac{1}{2}$. Below we
show that for any $H \in (0,1)$ using any function $\psi$ with two
vanishing moments in (\ref{eq:detail}) is sufficient for obtaining
an asymptotically Gaussian wavelet variance series
(\ref{eq:scale}).


\begin{theorem}\label{mainthm2}
Assume that $Z$ in (\ref{eq:detail}) is given by (\ref{logstock}),
with $\mu=0$ and $\sigma$ constant, and the analyzing wavelet in
(\ref{eq:detail}) has compact support and has  vanishing moments
of order at least 2. Then the logarithm of the scale spectrum
(\ref{eq:scale}), i.e $\log_{2}(S_{j})$, is asymptotically
normally distributed and satisfies the following asymptotic
relation (as $N \rightarrow \infty$):
\begin{equation}\label{nonuni}
\log_{2}S_{j} \sim \log_{2}(\sigma^2
K(H))+j(2H+1)+\frac{\epsilon_{j}}{\sqrt{N_{j}}\ln(2)} \qquad
j=1,...,\log_{2}(N)
\end{equation}
where $N_{j}=N/2^j$, $K(H)$ is given by (\ref{eq:kh}),
$\epsilon_{j}$ is $\mathcal{N}(0,1)$ and
\begin{equation}\label{eq:asycov}
cov\left(\frac{\epsilon_{j}}{\sqrt{N_{j}}},\frac{\epsilon_{i}}{\sqrt{N_{i}}}\right)\sim
\frac{1}{\sqrt{N_{j}N_{i}}},
\end{equation}
as $N \rightarrow \infty$.
\end{theorem}
First we will state a central limit theorem for heteroskedastic
random variables:
\begin{lemma}(\textbf{Berry-Essen Theorem}(see
\cite{kn:stroock})) Suppose $y_{1},y_{2},...,y_{n}$ are
independent random variables such that
\begin{alignat*}{2}
\EE\{y_{i}\} = 0, & \qquad \EE\{y_{i}^{2}\}= \sigma_{i}^{2} &
\qquad \text{and} \qquad \EE\{|y_{i}^{3}|\} =\rho_{i},
\end{alignat*}
and define
\begin{alignat*}{2}
s_{n}^{2}=\sum_{i=1}^{n}\sigma_{i}^{2} & \qquad \text{and} \qquad
r_{n}=\sum_{i=1}^{n}\rho_{i}.
\end{alignat*}
Let $F_{n}$ denote the distribution function of the normalized sum
$\sum_{i=1}^{n}y_{i}/s_{n}$. Then
\begin{equation*}
|F_{n}(x)-\Phi(x)|\leq 6 \frac{r_{n}}{s_{n}^{3}},
\end{equation*}
where $\Phi$ denotes the ${\cal N}(0,1)$ distribution function.
\end{lemma}
\emph{\textbf{Proof of Theorem \ref{mainthm2}:}}\newline
 First we will use the Berry-Essen
Theorem to show that $S_{j}$ given by ~\eqref{eq:scale} is
asymptotically normal ($N\rightarrow \infty$) with mean
proportional to $2^{j(2H+1)}$. Let $N_{j}=N/2^{j}$, and denote the
vector of scale coefficients at scale $j$ by
$d^{j}=[d_{j}(1),...,d_{j}(N_{j})]^{T}$. Also denote the
covariance matrix of $d^{j}$ by $C$. Note that $d^{j}$ has the
same law as $C^{1/2}\eta$, where $\eta$ is a vector of independent
${\cal N}(0,1)$ random variables. Let $A$ be the matrix that
diagonalizes $C$, {\em i.e.} $A^{T}CA=\Lambda$, where $\Lambda$ is
the matrix of eigenvalues of $C$. Since $\xi=A\eta$ has the same
distribution as $\eta$, we have
\begin{equation*}
N_{j}S_{j}={d^{j}}^{T}{d^{j}}\overset{d}{=}\eta^{T}C\eta
\overset{d}{=}\xi^{T}\Lambda\xi=\sum_{i}\lambda_{i}\xi_{i}^{2}.
\end{equation*}
On denoting $\tilde{\Lambda}=\Lambda/\EE\{S_{j}\}$, we have
\begin{equation*}
S_{j}=\EE\{S_{j}\}\bigg(1+ \frac{1}{N_{j}} \sum_{i=1}^{N_{j}}
\tilde{\lambda}_{i}(\xi_{i}^{2}-1)\bigg).
\end{equation*}
Define
\begin{equation*}
X_{j}=\frac{1}{N_{j}}\sum_{i=1}^{N_{j}}y_{i}
\end{equation*}
where $y_{i}=\tilde{\lambda}_{i}(\xi_{i}^{2}-1)$. The $y_{i}'s$
are independent random variables with the following properties:
\begin{equation*}
\EE\{y_{i}\}=0, \qquad \EE\{y_{i}^{2}\}=2\tilde{\lambda}_{i}^{2}
\qquad \text{end} \qquad \EE\{|y_{i}^{3}|\}\leq
28\tilde{\lambda}_{i}^{3},
\end{equation*}
which are easily derived from the fact that the moments of a
${\cal N}(0,\sigma^{2})$ random variable are given by the
following expression:
\begin{equation*}
\mu_{j}=\big(\frac{\sigma^{2}}{2}\big)^{j/2}
\frac{j!}{\frac{j}{2}!} , \qquad \text{for $j$ even},
\end{equation*}
and the odd moments are zero. By the Berry-Essen Theorem, it is
sufficient to show that
\begin{equation*}
J=\frac{\sum_{i=1}^{N_{j}}\tilde{\lambda}_{i}^{3}}{[\sum_{i=1}^{N_{j}}\tilde{\lambda}_{i}^{2}]^{3/2}},~\label{eq:J}
\end{equation*}
is small for large $N_{j}$. We first analyze the denominator.
\begin{equation*}
\sum_{i=1}^{N_{j}}\tilde{\lambda}_{i}^{2}=\frac{\sum_{n,m}(C_{nm})^{2}}{(\EE\{S_{j}\})^{2}}=\frac{N_{j}}{2}\bigg(\frac
{4}{N_{j}}\sum_{k=0}^{N_{j}-1}(N_{j}-k)[\rho(k)]^{2}-2\bigg)
\end{equation*}
where
\begin{equation*}
\rho(k)=\frac{\EE\{d_{j}(n)d_{j}(n-k)\}}{\EE\{[d_{j}(n)]^{2}\}}.
\end{equation*}
Let us introduce
\begin{equation*}
l(H)=\lim_{N_{j}\rightarrow \infty}\bigg(\frac
{4}{N_{j}}\sum_{k=0}^{N_{j}-1}(N_{j}-k)\rho(k)^{2}-2\bigg).
\end{equation*}
 We will now show that $\rho(k)$ decays as $k^{-2p+2H}$ (so that $l(H)$ is a constant), where $p$
is the number of vanishing moments and $H$ is the Hurst exponent
of the fBm.  We can write
\begin{equation*}
\EE\{d_{j}(n)d_{j}(n-k)\}=2^{-j}\sigma^{2}\EE\bigg\{\int dx \int
dt\, \psi_{j,n}(x) \psi_{j,n-k}(t) B^{H}(x) B^{H}(t)\bigg\}.
\end{equation*}
By ~\eqref{eq:covfbm} and Fubini's Theorem we have
\begin{equation*}
\EE\{d_{j}(n)d_{j}(n-k)\}=2^{-j-1}\sigma^{2}\int dx \int
dt\,\psi_{j,n}(x)
\psi_{j,n-k}(t)\big(|t|^{2H}+|x|^{2H}-2|t-x|^{2H}\big),
\end{equation*}
and since $\int \psi_{j,k}(x) \ dx=0$ for all $(k,j)$ then we have
\begin{equation}
\begin{split}
\EE\{d_{j}(n)d_{j}(n-k)\}&=-2^{-j}\sigma^{2}\int dx \int
dt\,\psi_{j,n}(x) \psi_{j,n-k}(t)|t-x|^{2H}\\
&=-2^{j(2H+1)}\sigma^{2}\int dt \int dx\, \psi(t)
\psi(x)|t-x+k|^{2H} ~\label{eq:reduced}.
\end{split}
\end{equation}
Using Taylor's formula we have,
\begin{equation}
\begin{split}
\bigg(1+\frac{t-x}{k}\bigg)^{2H}&=1+\sum_{q=1}^{2p-1} \frac{
\Gamma(2H+1)}{ \Gamma(2H-q+1)
\Gamma(q+1)}\bigg(\frac{t-x}{k}\bigg)^{q} \\
&+\frac{\Gamma(2H+1)}{\Gamma(2H-p+1)\Gamma(p+1)} \theta
(t,x,k)\bigg(\frac{t-x}{k}\bigg)^{2p}, ~\label{eq:gamma}
\end{split}
\end{equation}
where $\theta (t,x,k)<(1+a/k)$ and where $a$ is the support length
of the analyzing wavelet. Using \eqref{eq:reduced} and the facts
that the mother wavelet $\psi$ has vanishing
 moments of order $p$ and is compactly supported, we conclude that
 $\rho(k)$ decays as $k^{-2p+2H}$. (In \eqref{eq:gamma}, $\Gamma$
 denotes the gamma function.)
Therefore,
\begin{equation*}
l(H)=\lim_{N_{j} \rightarrow \infty}\bigg(\frac
{4}{N_{j}}\sum_{k=0}^{N_{j}-1}(N_{j}-k)\rho(k)^{2}-2\bigg)
\end{equation*}
is constant if $p\geq 2$. (Note that for any self-similar
stationary increment processes with self-similarity parameter $H$,
i.e. $X(at)\overset{d}{=}a^{H}X(t)$, for any $a>0$, $H$ must be in
(0,1]; see  \cite{kn:taqqu}.) Therefore having wavelets of
vanishing moments of order 2 is necessary to cover the range
$(0,1]$ for the Hurst parameter. Note that Haar (having $p=1$)
wavelets would work only for $H\in(0,\frac{3}{4}]$.

Now let us consider the numerator of ~\eqref{eq:J}. First we show
that the eigenvalues of $C/\EE\{S_{j}\}$ are bounded. By the
Gershgorin circle Theorem, the eigenvalue corresponding to a row
is not different from the corresponding diagonal element by more
than the sum of the other elements in the row, {\em i.e.},
\begin{equation}
|\tilde{\lambda_{i}}-1|\leq \sum_{n\neq i}|C_{in}|/C_{11}, \qquad
\text{where} \qquad C_{11}=\EE\{S_{j}^{2}\}. ~\label{eq:gers}
\end{equation}
Since $\rho(k)$ decays as $k^{2p-2H}$, the sum in ~\eqref{eq:gers}
approaches a constant in the limit as $ N_{j} \rightarrow \infty$
for any $H$ if the wavelet has at least two vanishing moments.
Therefore
\begin{equation*}
\lambda_{i} \leq K,
\end{equation*}
for some constant $K\geq 1$ independent of $N_{j}$. (Note that
$\lambda_{i}>0$.) Thus
\begin{equation*}
\sum_{i=1}^{N_{j}}\lambda_{i}^{3} \leq N_{j}K^{3}.
\end{equation*}
Hence, in the limit, $J$ in ~\eqref{eq:J} is given by
\begin{equation*}
J=\frac{(2K^{2}/l(H))^{3/2}}{\sqrt{N_{j}}},
\end{equation*}
and goes to zero. From the Berry-Essen Theorem we conclude that
\begin{equation*}
\frac{X_{j}}{\sqrt{l(H)}}=\frac{\sum_{i=1}^{N_{j}}\tilde{\lambda}_{i}(\xi_{i}^{2}-1)}{\sqrt{l(H)N_{j}}}
\end{equation*}
tends to a ${\cal N}(0,1)$ random variable in distribution. Thus,
asymptotically, $S_{j}$ is given by
\begin{equation*}
S_{j} \overset{d}{=} \EE\{S_{j}\}\bigg(1+\epsilon_{j}
\sqrt{\frac{l(H)}{N_{j}}}\bigg),
\end{equation*}
where $\epsilon_{j}$ is ${\cal N}(0,1)$.

One can also show that the $S_{j}$'s are asymptotically jointly
normal, by showing that $\sum_{j}a_{j}S_{j}$ is asymptotically
normal for any $(a_{j})_{1\leq j \leq M} $ (where $M$ is the
number of scales)
 using the same line of argument as above. From the
variance of the above sum one can find an expression for the
asymptotic normalized covariance of $(S_{j})$:
\begin{equation}
D_{j,i}=\frac{Cov(S_{j},S_{i})}{\EE\{S_{j}\}\EE\{S_{i}\}}
\underset{N\rightarrow \infty}{\sim} \frac{1}{\sqrt{N_{j}N_{i}}}
~\label{eq:covariance}
\end{equation}
where $N_{i}$ is the number of detail coefficients at scale $i$,
and $N$ is the total number of data points. The asymptotic
distribution of $\log_{2}(S_{j})$ can be derived exactly the same
way as in \cite{kn:solna} for the pure fBm case, with Haar
wavelets as the analyzing wavelets; therefore we will not repeat
this analysis here. The distribution of $\log_{2}(S_{j})$ is given
by
\begin{equation*}
\log_{2}(S_{j})\overset{d}{=}\log_{2}(\EE\{S_{j}\})+\epsilon_{j}\sqrt{\frac{l(H)}{N_{j}\ln(2)}},
\end{equation*}
where $\epsilon_{j}$ is ${\cal N}(0,1)$. \hfill $\Box$

In view of Theorem (\ref{mainthm2}) we can use the generalized
least squares estimate to estimate the Hurst parameter. If we
denote $c:=log_{2}(\sigma^2 K(H))$, and $h:=2H+1$ then the
generalized least squares estimate of $b=[c,h]^{T}$ is given by
\begin{equation}\label{eq:esteq}
 \tilde{b}=(X^{T}D^{-1}X)^{-1}X^{T}D^{-1}M
\end{equation}
where $M=[\log_{2}(S_{1}),...,\log_{2}(S_{J})]^{T}$,
($J=\log_{2}(N)$), $D$ is given by ~\eqref{eq:covariance}, and $X$
is given by
\begin{equation*}
X=
\begin{bmatrix}
1 & 1 \\ 1 & 2 \\ \vdots & \vdots \\ 1 & J.
\end{bmatrix}
\end{equation*}

We have
\begin{equation*}
\begin{split}
\EE\{\tilde{c}\}&= \log_{2}(\sigma^{2}K(H)),\\
\EE\{\tilde{h}\}&=2H+1,
\end{split}
\end{equation*}
and
\[
  \EE\{\tilde{b}\tilde{b}^{T}\}=\left(X^{T}D^{-1}X
  \ln(2)^{2}\right)^{-1} \sim \frac{1}{N}
\]
for large $N$.  In view of these we have the following estimator
for the Hurst parameter:
\begin{equation}\label{finally}
\tilde{H}=(\tilde{h}-1)/2,
\end{equation}
with variance
\begin{equation}
Var(\tilde{H}) \sim 1/(4N)
\end{equation}
for large $N$.

In the next section we will allow $\sigma$ in (\ref{logstock}) to
be stochastic. In particular we take $\sigma$ to be any stochastic
process having sufficient regularity.
\subsection{Robustness to Seasonalities and Volatility
Persistence}\label{presistence}

We first present an empirical verification of robustness to
seasonalities which is followed by a theoretical verification of
robustness both to seasonalities and to volatility persistence.

\begin{figure}
\begin{center}
\includegraphics[width = 0.8\textwidth,height=5cm]{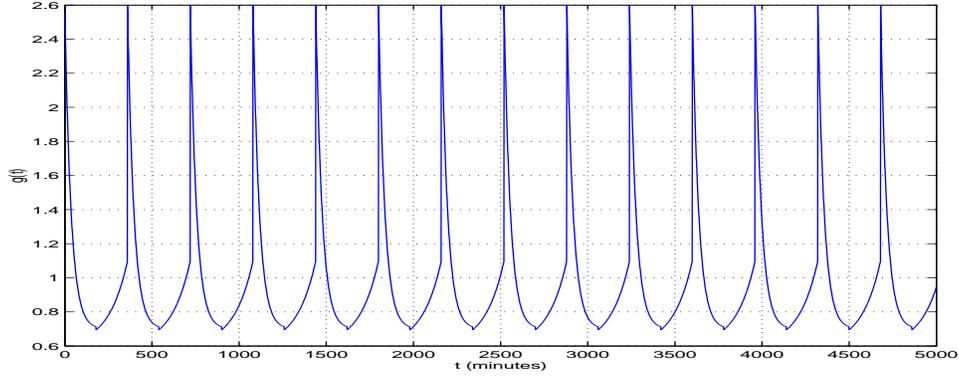}
\caption{The periodic function $g$ estimated from the S\&P 500,
whose one period is a good representative of the intraday
variability of the index.} \label{fig:periodic}
\end{center}
\end{figure}
We will denote the seasonal component with $g(x)$, where $x>0$ is
the time from the beginning of the segment under discussion.
(Since the seasonal component is deterministic, we will denote it
by $g(t)$ instead of $\sigma_{t}$ to avoid confusion.) In examples
such as $g(x)=(x+b)^q$, $b$ denotes the beginning of the segment.

\begin{figure}
\begin{center}
\includegraphics[width = 0.8\textwidth,height=5cm]{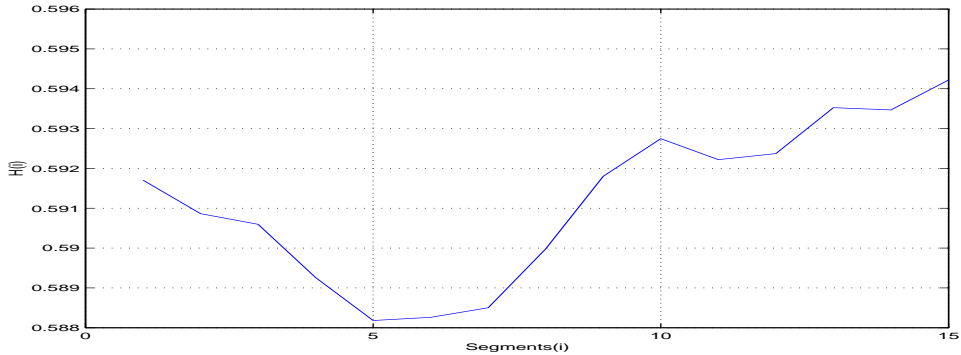}
\caption{Hurst estimates from simulated data for ~\eqref{eq:model}
with $g(x)$ from Fig.~\ref{fig:periodic}, and with segments of
length $2^{15}$.} \label{fig:periodich}
\end{center}
\end{figure}

When (\ref{eq:model}) is implemented with a periodic $g$ given by
Fig.~\ref{fig:periodic}, which represents actual intraday
variability, then the Hurst estimates do not change significantly.
(The intraday variability envelope $g$ of Fig. \ref{fig:periodic}
was estimated from the S\&P 500 index data as in
\cite{kn:fouque}.) The Hurst estimates also do not change when the
amplitude of this periodic component is changed be a factor of
100.  Estimates of the Hurst parameter for segments of length
$2^{15}$ are shown in Fig.~\ref{fig:periodich}. (Compare with Fig.
\ref{fbmhurts}). The mean of the estimates over the segments is
0.5912, and the standard deviation is 0.0020.

\begin{figure}
\begin{center}
\includegraphics[width = 0.8\textwidth,height=5cm]{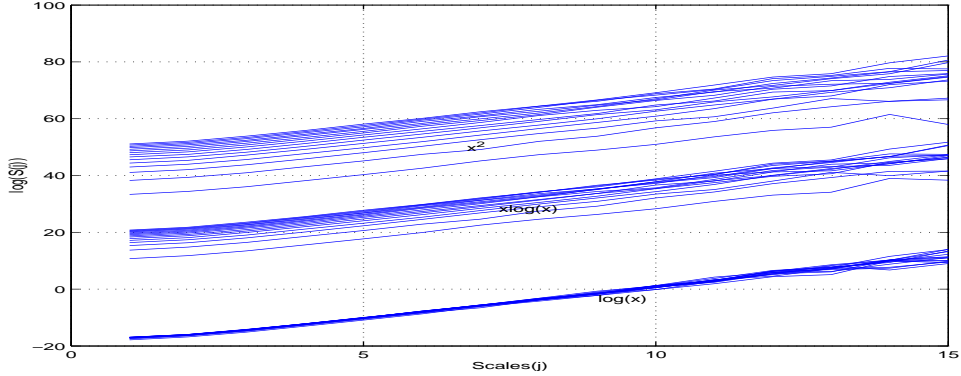}
\caption{The log-scale spectrum for ~\eqref{eq:model} with
$g(x)=\log(x+b)$, $g(x)=(x+b)\log(x+b)$ and $g(x)=(x+b)^{2}$.}
\label{fig:scaleforvarious}
\end{center}
\end{figure}

The scale spectra of (~\eqref{eq:model}) corresponding to various
$g$'s, namely $g(x)=\log(x+b)$, $g(x)=(x+b)\log(x+b)$,
$g(x)=(x+b)^{2}$, are plotted together on the same graph for
comparison (Fig.~\ref{fig:scaleforvarious}). One immediately
notices that the slope of the scale spectrum is invariant to the
choice of $g$ for these examples; only the amplitude changes with
$g$. For the realization of fBm shown in Fig.~\ref{fig:fbm}, the
mean and the std. for this parameter (over the segments) are:
0.5930, 0.0088; 0.5942, 0.0642; and 0.5953, 0.0062 for
$g(x)=\log(x+b)$, $g(x)=(x+b)\log(x+b)$, and $g(x)=(x+b)^{2}$
respectively.

Path properties play an important role in the robustness of the
estimator developed in the previous section in the case of
stochastic volatility. First note that paths of $B^{H}$ are almost
surely H\"{o}lder continuous of order $\lambda$ for all
$\lambda<H$, due to the Kolmogorov-\u{C}entsov Theorem
(\cite{kn:karat}).

The following result due to Ruzmaikina \cite{kn:ruzmaikina} and
Z\"{a}hle \cite{kn:zahle} gives the path properties of the
stochastic integrals of a certain class of processes with respect
to fBm.
\begin{lemma}\label{integral}
Suppose $\sigma$ is a stochastic process with almost surely
H\"{o}lder continuous paths of order $\gamma>1-H$ on the interval
$[0,T]$. Then the integral
\begin{equation}\label{fBmint}
I_{t}=\int_{0}^{t}\sigma_{s}dB^{H}_{s}
\end{equation}
exists almost surely as a limit of Riemann-Stieltjes sums.
Furthermore, the process I is almost surely $\beta$-H\"{o}lder
continuous on $[0,T]$ for any $\beta<H$.
\end{lemma}

Note that $\sigma$ does not have to be adapted with respect to the
natural filtration of $B^{H}$. For $H>1/2$ an example of $\sigma$
satisfying the conditions of Lemma \ref{integral} is the Wiener
process. Any continuous periodic function also satisfies the
assumptions of this theorem. (This is a rather straightforward
example, but we cite it due to its relevance to seasonality
issue.)

The following lemma is the key result for robustness; it gives a
bound on the wavelet detail coefficients (\ref{eq:detail}) for
functions with certain regularity.

\begin{lemma}(See \cite{kn:mallat})\label{jaffard}
A function $f$ is H\"{o}lder continuous of order $\gamma$ if and
only if the scale coefficients corresponding to $f$ satisfy
\begin{equation*}
|d_{j}(k)| \leq A 2^{j(\gamma+1/2)}, \quad \forall \ (j,k) \in
\mathbb{Z}^{2},
\end{equation*}
for some $A<\infty$.
\end{lemma}

The H\"{o}lder continuity exponent of a function is related to its
finer scales; therefore one must use the scale coefficients
defined in (\ref{eq:detail}) for $j<0$. Using Lemma \ref{jaffard}
we have the following result.

\begin{lemma}
Suppose f is a function that is H\"{o}lder continuous of order
$\lambda<H$ and $S_{j}$ is its scale spectrum. Then
\begin{equation}\label{inequality}
S_{j} \leq C_{\gamma}2^{j(2\gamma+1)}, \ \ \forall j \in
\mathbb{Z},\ \forall \gamma<H,
\end{equation}
for some $C_{\gamma} \in (0,\infty)$ and moreover
(\ref{inequality}) does not hold for $\gamma \geq H$ for
infinitely many $j \in \mathbb{Z}_{-}$
\end{lemma}
and
\begin{equation*}
\liminf_{j \rightarrow -\infty}\frac{\log(S_{j})}{j}=2H+1.
\end{equation*}

Let us summarize the results of this section in the following
theorem.

\begin{theorem}\label{mainthm1}
Suppose $\sigma$ is a stochastic process with almost surely
H\"{o}lder continuous paths of order $\gamma>1-H$ on the interval
$[0,T]$. Then there exists a random variable  $C_{\gamma}(\omega)
\in (0,\infty)$ such that the scale spectrum of the integral
(\ref{fBmint}) satisfies
\begin{equation}\label{inequality1}
S_{j}(\omega) \leq C_{\gamma}(\omega)2^{j(2\gamma+1)}, \ \ \forall
j \in \mathbb{Z},\ \forall \gamma<H,
\end{equation}
almost surely. Moreover (\ref{inequality1}) almost surely does not
hold for $\gamma > H$ for infinitely many $j \in \mathbb{Z}_{-}$
and
\begin{equation*}
\liminf_{j \rightarrow -\infty}\frac{\log(S_{j}(\omega))}{j}=2H+1,
\end{equation*}
almost surely.

\end{theorem}




In Section \ref{section4}, the domain of the wavelet is taken to
be on the order of the mesh size of the discrete samples of the
data. However the sample path properties show themselves in the
finer detail coefficients ($j<0$ in (\ref{eq:detail})). Then,
using the fact that for any function $f$ that is
$\lambda$-H\"{o}lder continuous, $g(t)=f(at)$ is also
$\lambda$-H\"{o}lder continuous, it can be seen that a scale
spectrum for the finer scales can be obtained from a scale
spectrum corresponding to coarser scales. Letting
$J=\log_{2}(N)+1$, where $N$ is as in (\ref{eq:scale}), we define
$\psi'(t)=2^{-J/2}\psi(at/2^{J})$. Then $\psi'$ is a wavelet which
has the same number of vanishing moments as $\psi$. Taking $a=1$,
and defining the scale coefficients as
\begin{equation*}
d'_{j}(k):=2^{-j/2}\int_{-\infty}^{+\infty}\psi'(2^{-j}t-k)Z(t)
dt, \quad j \in \mathbb{Z}_{-}
\end{equation*}
and scale spectrum as
\begin{equation*}
S'_{j}:=\frac{1}{N/2^{j}}\sum_{k=1}^{N/2^{j}}[d'_{j}(k)]^{2} \quad
j \in \mathbb{Z}_{-}
\end{equation*}
we have $S'_{j}=S_{J+j}$. So as the number of samples of the data
increases ($J \rightarrow \infty$), we can consider finer and
finer detail coefficients and the corresponding scale spectrum.

Since the log scale spectra corresponding to the synthetically
constructed data for various kinds of seasonalities (see Fig.
\ref{fig:scaleforvarious}) and corresponding to the real data were
linear, from Theorem \ref{mainthm1}
we can conclude that a linear regression is an accurate way of
estimating the H\"{o}lder continuity exponent $H$ of a sample
path. This technique has been employed by Arneodo
\cite{kn:arneodo} for estimating the multifractal spectrum of a
given sample path. (Also see \cite{kn:mallat}). Since the
estimator given by (\ref{eq:esteq}) and (\ref{finally}) is a
linear weighted least squares fit to the scale spectrum (giving
more weight to the smaller scales, the weighting factor being
proportional to the number of scale coefficients at the given
scale scale) it is equal to the multifractal spectrum estimator of
Arneodo.

\section{Hurst Parameter Estimation for the S\&P 500}
After segmenting the samples of $\log(P_{t}/P_{0})$ from our S\&P
500 data set into dyadic segments, we estimated the Hurst
parameter for each of the segments using the estimator given by
(\ref{finally}). (Note that $N$ is the number of points in a given
segment.)

For each segment (\ref{finally}) requires computing the scale
spectrum (\ref{eq:scale}) which further requires the computation
of the detail coefficients (\ref{eq:detail}) for every scale. It
may seem that this is computationally expensive for high frequency
data; however due to the pyramidal algorithm described in Section
\ref{pyramidalalgorithm}, this is not an issue. The pyramidal
algorithm calculates the wavelet coefficients for any number of
scales using octave filter banks given the initial approximation
coefficients. Therefore the detail coefficients need only be
computed at the initial scale. The detail coefficients at higher
scales are computed from these initial coefficients via the
pyramidal algorithm, which uses only the approximate coefficients
of the preceding scale for calculating the detail coefficients of
the next scale.

 In our model we want to introduce the flexibility of having
a variation in $H$. If we further partition the segments into
smaller segments of equal length and estimate the Hurst parameter
for these smaller segments, we expect to see noise in our
estimates due to the noise introduced by making the segmentation
length smaller. To be able to make a comparison between the Hurst
estimates corresponding to different segmentation levels we must
filter out the extra noise introduced. The following section
introduces a method to remove this finite segmentation noise.

\subsection{Filtering the Finite Segmentation Noise}
We know that the log-scale spectrum method yields an
asymptotically efficient estimator, and thus having smaller
segment lengths will introduce noise into the estimates. To deal
with the noise due to finite segmentation length, we will follow
the approach of Papanicolaou and Solna \cite{kn:solna}, which we
now review.

On letting $\hat{h}^{i}$ denote the slope of the log-scale
spectrum for the $i$ th segment, we model it as
\begin{equation*}
\hat{h}^{i}=h^{i}+\zeta^{i},
\end{equation*}
where $h^{i}$ is the true slope for the $i$th segment
$(h^{i}=2H^{i}+1)$, and $\zeta^{i}$ is a random variable that
models the finiteness of the segments. The error term will be
assumed to have zero mean, and for different segments the error
terms will be assumed to be uncorrelated, i.e.
\begin{equation*}
\EE\{\zeta^{i}\zeta^{j}\}=0, \qquad i\neq j.
\end{equation*}
Here it is assumed that $(H^{i})$ is a stationary stochastic
process independent of the fBm. On assuming that the slope process
is exponentially correlated, $l_{h}$ denotes its correlation
length, and $\sigma_{h}^{2}$ denotes its variance, we have
\begin{equation*}
C_{h}(i,j)=\EE\big\{(h^{i}-\EE\{h^{i}\})(h^{j}-\EE\{h^{j}\})\big\}=\sigma_{h}^{2}\exp(-L|i-j|/l_{h})~\label{eq:slopecov}
\end{equation*}
where $L$ is the segment length.

Here we will give a minimum variance unbiased linear estimator for
the slope process $(\hat{h}^{i})$. Let $\hat{K}$ denote the vector
whose components are the estimates $(\hat{h}^{i})$, and let $K$
denote the vector whose components are the corresponding
realizations of the slope process $(h^{i})$. We want to find a
filter such that
\begin{equation*}
\EE\{\parallel \Gamma\hat{K}-K \parallel^{2}\}
\end{equation*}
is minimized under the constraint that the mean is preserved,
i.e.,
\begin{equation*}
\Gamma \EE\{K\}=\EE\{K\}.
\end{equation*}
It can be shown that $\Gamma$ is given by
\begin{equation}
\Gamma=(C_{h}+C_{\zeta})^{-1}\big[C_{h}+u^{T}\otimes
\EE\{K\}\big], ~\label{eq:gamma2}
\end{equation}
where the vector $u=(u_{i})$ is given by
\begin{equation*}
u_{i}=
\frac{\EE\{h^{i}\}-\EE\{K^{T}\}(C_{h}+C_{\zeta})^{-1}C_{h,i}}{\EE\{K^{T}\}(C_{h}+C_{\zeta})^{-1}\EE\{K\}}
\end{equation*}
and $C_{\zeta}$ is the diagonal covariance matrix of the
estimation errors $\zeta^{i}$. Here $C_{h,i}$ denotes the $i$th
column of $C_{h}$. To be able to implement $\Gamma$ one must
estimate $\sigma_{h}^{2}$ and $l_{h}$ of ~\eqref{eq:slopecov}, and
the variance $\sigma_{w}^{2}$ of the noise process. For this
purpose we will examine the empirical variogram of the slope
process $(\hat{h}^{i})$. The variogram at lag $j$ is given by
\begin{equation*}
V(j)=\frac{1}{2(J-j)}\sum_{k=1}^{J-j}(\hat{h}^{k+j}-\hat{h}^{k})^{2}
\end{equation*}
where $J$ is the number of segments. Since
\begin{equation}
\EE\{V(j)\}=\sigma_{h}^{2}(1-\exp(-L|j|/l_{h})+\sigma_{\zeta}^{2},
~\label{eq:vario}
\end{equation}
fitting (by a weighted least squares fit) the left-hand side of
~\eqref{eq:vario} to the empirical variogram yields the estimates
for $\sigma_{h}$, $\sigma_{\zeta}$ and $l_{h}$. Here it should be
noted that the initial values must be chosen carefully. The most
important parameter seems to be the mean correlation length
$l_{h}$, and we choose it to be longer than the segments used in
estimating the slopes in order to have approximate stationarity
relative to segmentation. Also note that it is necessary to
perform a weighted fit to the empirical variogram, because there
are finitely many segments, and therefore the empirical variogram
is closer to its expected value for smaller lags.

\begin{figure}
\begin{center}
\includegraphics[width = 0.8\textwidth,height=6cm]{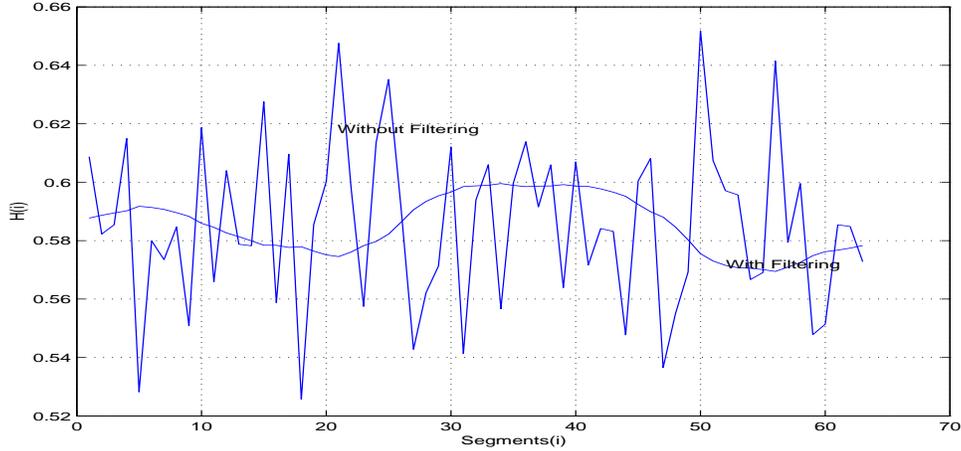}
\caption{Hurst estimates for the realization given in
Fig.~\ref{fig:fbm} for 64 segments of length $2^{13}$ with and
without filtering. } \label{fig:comp}
\end{center}
\end{figure}

We will now illustrate the power of filtering in removing the
effects of noise due to finite segmentation length on one
realization of the synthetically created fBm of
Fig.~\ref{fig:fbm}. We partitioned the data into segments of
length $2^{13}$, estimated the Hurst parameter for each of the
segments, and then applied filtering. ($\Gamma P$ where $\Gamma$
is given in ~\eqref{eq:gamma2}.) The Hurst estimates we obtained
with and without filtering are given in Fig.~\ref{fig:comp}. The
standard variation without filtering is 0.0283, and the standard
variation after filtering is 0.0097; so clearly we mitigate to the
finite segment length effects by filtering.

\subsection{Results on the S\&P 500 Index}
We now turn our attention to the analysis of the S\&P 500 index.
As noted above, we consider data taken at one-minute intervals
  over the
course of 11.5 years from January 1989 to May 2000. (We take the
closing price of each minute.) The data consists of 1,128,360
observations, which is on the order of $2^{20}$.

\begin{figure}
\begin{center}
\includegraphics[width = 0.8\textwidth,height=6cm]{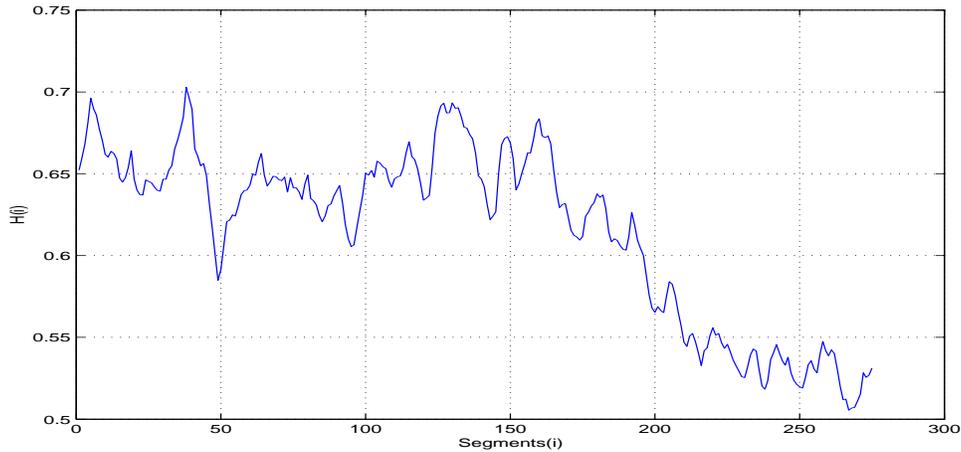}
\caption{Hurst parameter estimates for the S\&P 500 data with
 segment lengths of
$2^{12}$} \label{fig:h12}
\end{center}
\end{figure}

When the data is segmented into 275 segments of length $2^{12}$
(approximately two weeks) and the above methodology is applied,
 we obtain the Hurst parameter estimates
shown in Fig.~\ref{fig:h12}. (The mean is 0.6156, and the standard
deviation is 0.0531, which supports the idea of local variation,
i.e. the Hurst parameter varies significantly from segment to
segment.)

\begin{figure}
\begin{center}
\includegraphics[width = 0.8\textwidth,height=6cm]{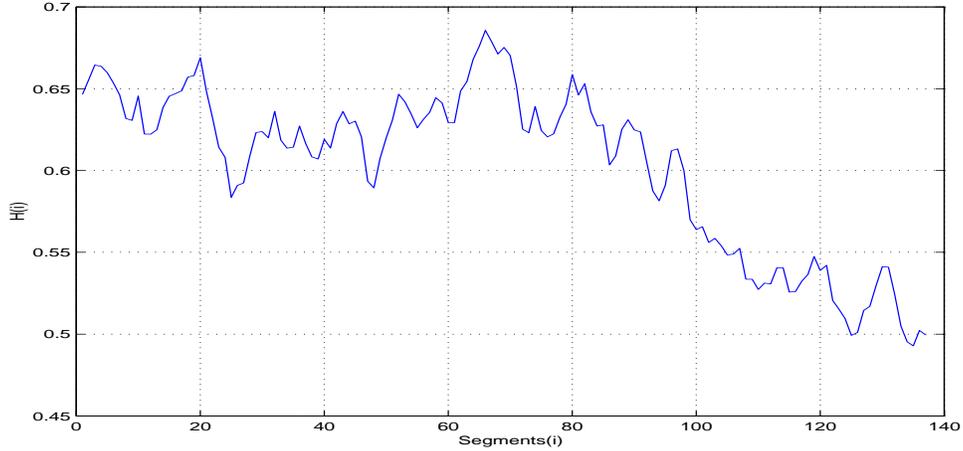}
\caption{Hurst parameter estimates for the S\&P 500 data with
 segment lengths of
$2^{13}$} \label{fig:h13}
\end{center}
\end{figure}

Alternatively, when the data is segmented into 137 segments of
length $2^{13}$ (approximately four weeks), we obtain the Hurst
parameter estimates shown in Fig.~\ref{fig:h13}. (The mean is
0.6027, and the std. is 0.0504.)

\begin{figure}
\begin{center}
\includegraphics[width = 0.8\textwidth,height=6cm]{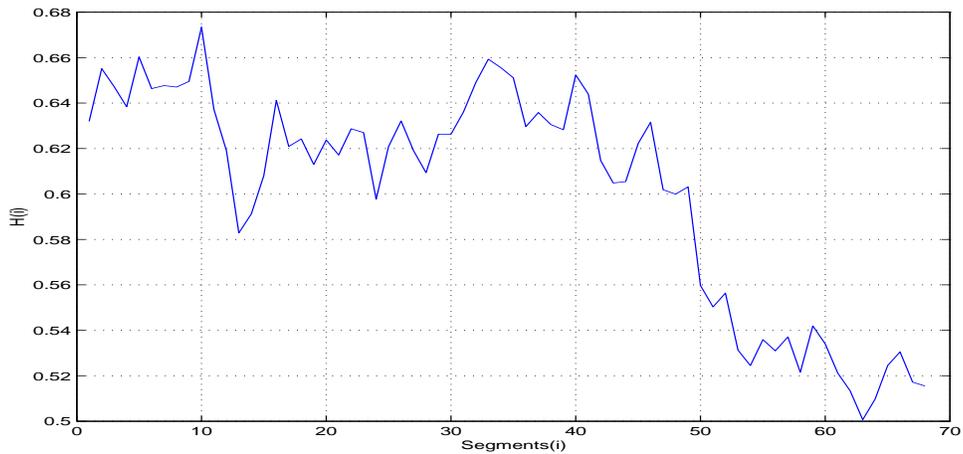}
\caption{Hurst parameter estimates for the S\&P 500 data with
 segment lengths of
$2^{14}$} \label{fig:h14}
\end{center}
\end{figure}

\begin{figure}
\begin{center}
\includegraphics[width = 0.8\textwidth,height=6cm]{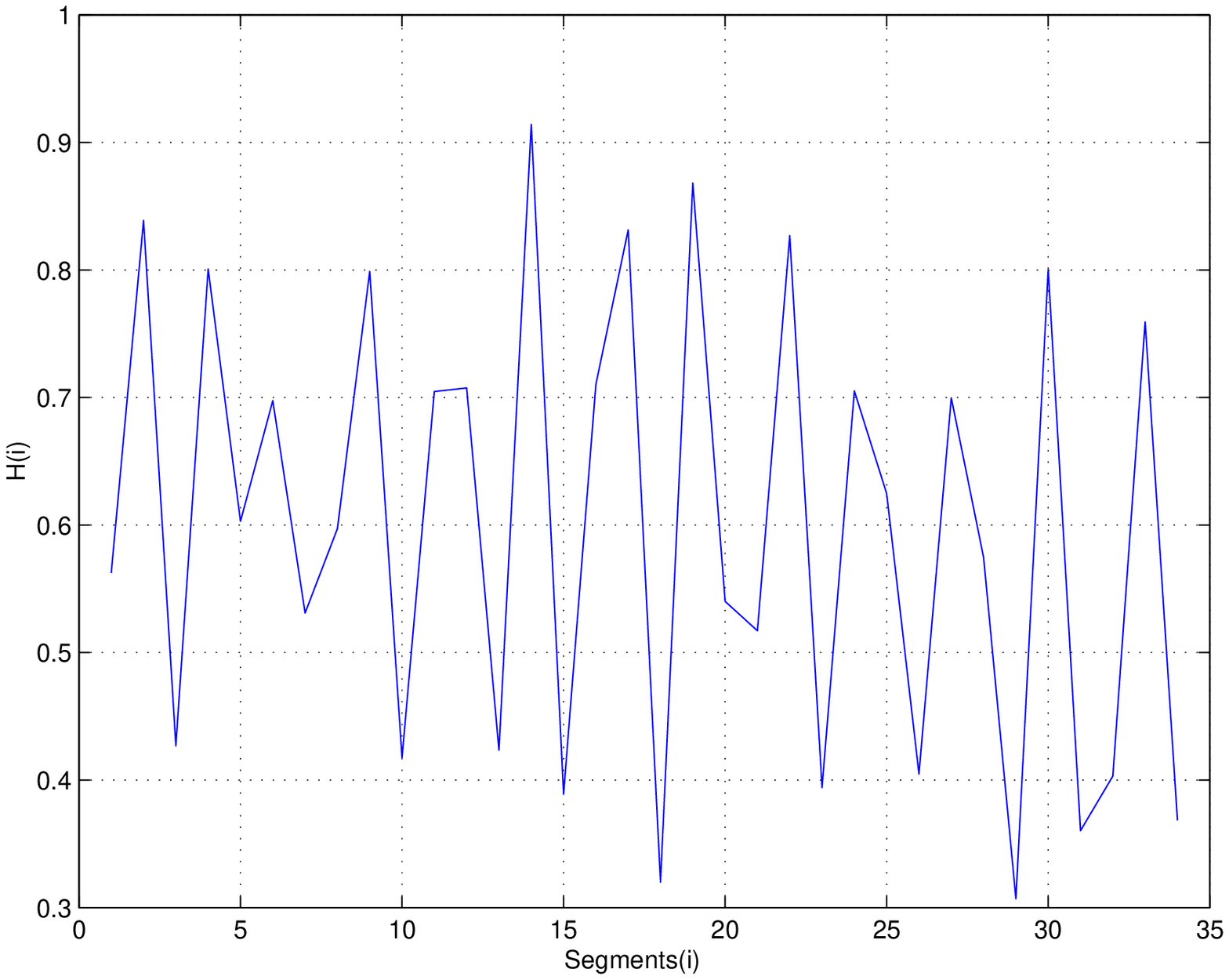}
\caption{Hurst parameter estimates for the S\&P 500 data with
 segment lengths of
$2^{15}$} \label{fig:h15}
\end{center}
\end{figure}

Similarly, when the data is segmented into 68 segments of length
$2^{14}$ (approximately four weeks), we obtain the Hurst parameter
estimates shown in Fig.~\ref{fig:h14}. (The mean is 0.6011 and the
std. is 0.0487.) And,  finally, when the data is segmented into 34
segments of length $2^{15}$ (approximately eight weeks), we obtain
the Hurst parameter estimates given in Fig.~\ref{fig:h15}. (The
mean is 0.6008 and the std. is 0.1821.)

\begin{figure}
\begin{center}
\includegraphics[width = 0.8\textwidth,height=6cm]{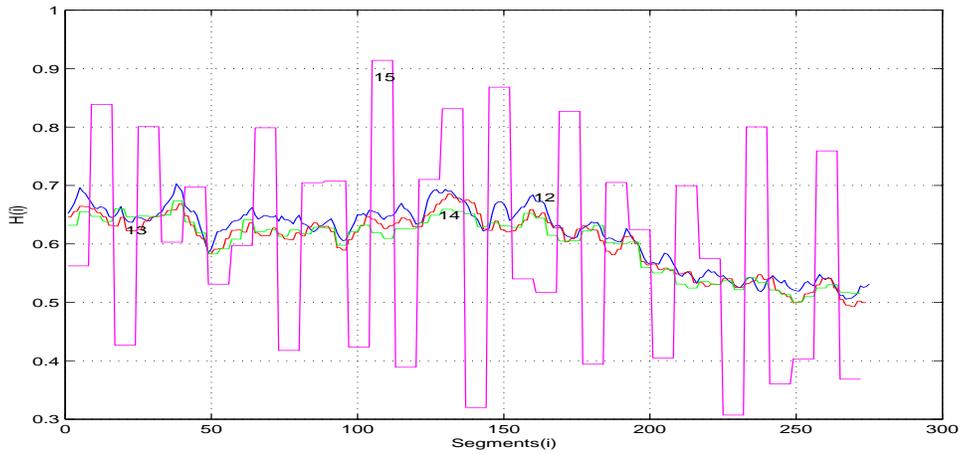}
\caption{Hurst parameter estimates for the S\&P 500 data with
 various segment lengths.
The strongly varying graph corresponds to estimates for segment
lengths of $2^{15}$ points. The graphs corresponding to the
segment lengths of $2^{12}$, $2^{13}$, $2^{14}$ can be
distinguished from the line intensity, where the intensity
decreases as the segment length increases, moreover the graphs are
labeled by $log_{2}$ of the corresponding segments length. }
\label{fig:hcomparison}
\end{center}
\end{figure}

If we plot these results on the same axes as shown in
Fig.~\ref{fig:hcomparison} we will arrive at the significant
observation that the length of a stationary segment is $2^{14}$,
which corresponds to approximately 2 months. That is, when the
segments are of length $2^{15}$, the nonstationarity is dominant.

\subsection{Increase in the Market Efficiency}

From Fig.~\ref{fig:h14}, one sees that, although the Hurst
parameter of this data set is significantly above the efficient
markets value of $H=\frac{1}{2}$, it began to approach that level
in 1997 (segment $50$). We conjecture that this behavior of the
market might be related to the increase in Internet trading, which
had the three-fold effect of increasing the number of small
traders, the frequency of trading activity and the availability of
market data. This observation is modeled in \cite{kn:bhs}, with a
simple microstructure model for the price evolution of a financial
asset where the price is driven by the demand of many small
investors whose trading behavior exhibits ``inertia''. This means
that the agents trade the asset infrequently and are inactive most
of the time. It is shown that when the price process is driven by
market imbalance, the logarithm of the price process is
approximated by a process of the form (\ref{logstock}). Moreover
it is shown that as the frequency of trading increases, the price
process can be approximated by geometric Brownian motion, which is
consistent with the above comments.



\section{Conclusion}
In this paper we have developed a method to investigate long range
dependence, which is quantified by the Hurst parameter, in high
frequency financial time series. Our method  exhibits robustness
to the non-stationarities that are present in the data, e.g.
seasonal volatility, fat-tailed distributions of the increments,
and possible variations in the Hurst parameter. (In fact, the
Hurst parameter reflects the relative frequency of the trading
activity of the market participants, and hence variations in the
Hurst parameter are expected \cite{kn:bhs}.) The segments of
stationarity for the Hurst parameter are byproducts of this
analysis. They are found to be approximately two months in
duration for S\&P 500 index data sampled at one minute intervals.
Strikingly, the Hurst parameter was around the $0.6$ level for
most of the 1990s, but dropped closer to the efficient markets
level of $0.5$ in the period 1997-2000, coinciding with the growth
in Internet trading among small investors.


\setcounter{section}{0}%
\renewcommand{\thesection}{\Alph{section}}%
\section{Appendix}
\subsection{Wavelet Theory}\label{wavelet}
For a more detailed treatment see ~\cite{kn:icd},
~\cite{kn:mallat} or ~\cite{kn:strang}.
\subsubsection{Multi-resolution Analysis}\label{multires} A wavelet $\psi$ is a
function mapping $\IR$ to $\IR$ such that the dilated and
translated family
\begin{equation*}
\psi_{j,k}(t)=2^{-j/2}\psi(2^{-j}t-k) \qquad \text{for} \ \
(j,k)\in\mathbb{Z}^{2}, ~\label{eq:psi}
\end{equation*}
is an orthonormal basis of $L^{2}(\IR)$. A wavelet is defined via
a \emph{scaling function} through multi-resolution analysis (MRA).
A sequence of closed subspaces $\{V_j, j\in\mathbb{Z}\}$ of
$L^{2}(\IR)$ is an MRA if the following six properties are
satisfied:
\begin{gather*}
 \bigcap_{j \in \mathbb{Z}}V_{j}=\{0\},\\
 \mbox{Closure}(\bigcup_{j \in \mathbb{Z}}V_{j})=L^{2}(\IR),\\
V_{j+1}\subset V_{j}, \\\forall(j,k) \in \mathbb{Z}^{2}, \ \  f(t)
\in V_{j} \Leftrightarrow f(t-2^{j}k) \in V_{j}, \\ \forall j \in
\mathbb{Z}, \ \ f(t) \in V_{j} \Leftrightarrow f(\frac{t}{2}) \in
V_{j+1},
\end{gather*}
and, there exists a function $\phi(t)$ in $V_{0}$, called the
scaling function, such that the collection
\[ \{\phi(t-k),k \in
\mathbb{Z}\} \] is a Riesz basis for $V_{0}$.

It follows that the scaled and shifted functions of the scaling
function  \[\{\phi_{j,k}(t)=2^{-j/2}\phi(2^{-j}t-k),k \in
\mathbb{Z}\}\] is an orthonormal basis of $V_{j}$ for all $j$.

Orthonormal wavelets carry the details necessary to increase the
resolution of a signal approximation. The approximations of a
function $f\in L^{2}(\IR)$ at scales $2^{j}$ and $2^{j-1}$ are
respectively equal to its orthogonal projections onto $V_{j}$ and
$V_{j-1}$. Let $W_{j}$ be the orthogonal complement of $V_{j}$ in
$V_{j-1}$, {\em i.e.} $V_{j-1}=V_{j}\oplus W_{j}$. Then the
orthogonal projection of $f$ onto $V_{j-1}$ can be decomposed as
the sum of orthogonal projections onto $V_{j}$ and $W_{j}$. The
projection onto $W_{j}$ provides the \emph{details} that appear at
scale $2^{j-1}$ but which disappear at the coarser scale $2^{j}$.
One can construct an orthonormal basis of $W_{j}$ by scaling and
translating a wavelet $\psi \in V_{0}$, and show that the family
given by ~\eqref{eq:psi} is an orthonormal basis for $L^{2}(\IR)$.
Since $\psi$ and $\phi$ are in $V_{0}$, and by the properties of
MRA we have:
\begin{align*}
\phi(t)&=\sqrt{2}\sum_{k=0}^{M}h(k)\phi(2t-k),\\ \intertext{and}
\psi(t)&=\sqrt{2}\sum_{k=0}^{M}g(k)\phi(2t-k),
\end{align*}
where $h$ is a low-pass filter satisfying some admissibility
conditions (see  \cite{kn:mallat}), and $g$ is the conjugate
mirror filter of $h$:
\begin{equation*}
g(n)=(-1)^{1-n}h(1-n).
\end{equation*}
Therefore wavelets are specified via the scaling filter $h$.

Wavelets are capable of removing nonstationarities because they
have vanishing moments. We say that $\psi$ has $p$ vanishing
moments if it is orthogonal to any polynomial of degree less than
$p$:
\begin{equation}\label{vmoments}
\int_{-\infty}^{+\infty}t^{k}\psi(t)dt=0 \qquad \text{for} \qquad
0\leq k <p.
\end{equation}
The most versatile wavelet family is the family of Daubechies
compactly supported wavelets, which are enumerated by their number
of vanishing moments. The Daubechies compactly supported wavelet
with $p=1$ is the Haar wavelet, which is the only wavelet in this
family for which an explicit expression can be found. In our
analysis we used Daubechies compactly supported wavelets with
$p=2$.

\subsubsection{Pyramidal Algorithm (Mallat Algorithm) and its
Initialization}\label{pyramidalalgorithm} Let us denote the
projection of a function $f\in L^{2}(\IR)$ onto $V_{j}$ and
$W_{j}$ respectively by
\begin{equation*}
a_{j}(k)=<f,\phi_{j,k}> \ \ \text{and} \ \
d_{j}(k)=<f,\psi_{j,k}>, \quad k \in \mathbb{Z},
\end{equation*}
where $<\cdot,\cdot>$ denotes the standard $L^2$ inner product.
The pyramidal algorithm calculates these coefficients efficiently
with a cascade of discrete convolutions and subsamplings. Denote
time reversal by $\bar{x}(n)=x(-n)$ and upsampling by
\begin{equation*}
\check{x}(n)=\begin{cases} x(n)& \text{if $n$ is even}, \\ 0&
\text{if n is odd}.
\end{cases}
\end{equation*}
The pyramidal algorithm is then given by the following theorem:
\begin{theorem}(Mallat  \cite{kn:mallat})
Decomposition:
\begin{align*}
a_{j+1}(p)&=\sum_{n=-\infty}^{+\infty}h(n-2p)a_{j}(n)=a_{j}\star
\bar{h}(2p) \\
d_{j+1}(p)&=\sum_{n=-\infty}^{+\infty}g(n-2p)a_{j}(n)=a_{j}\star
\bar{g}(2p).
\end{align*}
Reconstruction:
\begin{equation*}
\begin{split}
a_{j}(p) & =
\sum_{n=-\infty}^{+\infty}h(p-2n)a_{j+1}(n)+\sum_{n=-\infty}^{+\infty}g(p-2n)d_{j+1}(n)
\\
&=\check{a}_{j+1} \star h(p)+\check{d}_{j+1} \star g(p).
\end{split}
\end{equation*}
\end{theorem}
To compute the detail coefficient at scale $j$, we use only the
approximation coefficient at the previous scale $a_{j-1}$.
 Note that the domain of $h$ and $g$ are compact if
we use Daubechies wavelets with compact support.

 The pyramidal algorithm assumes initially that it is given the
wavelet coefficients at a fine scale, and proceeds to compute the
detail coefficients at higher scales. The initial sequence
$a_{0}(k)$ requires the evaluation of a continuous time integral,
\begin{equation}
a_{0}(k)=\int_{\IR}f(t)\phi(t-k)\,dt,~\label{eq:int}
\end{equation}
where $\phi$ is the scaling function.

 Typically what is done is to set $a_{0}(k)=f(k)$, an ad-hoc
procedure that will almost certainly introduce errors. (An
exception is the case in which \emph{coiflets} \cite{kn:mallat}
are used, since in that case the scaling function has vanishing
moments of the same order as the wavelet.) Here, we will replace
the continuous time integral by a sum:
\begin{equation*}
a_{0}(k)=\sum_{n}f(n)\phi(n-k).
\end{equation*}
Note that this sum is equal to the integral of ~\eqref{eq:int}
when $f(t)$ is a low order polynomial. (See  \cite{kn:strang}.) An
explicit expression for $\phi$ is not known, however $\phi$ at the
integer points can be calculated from the defining recursion for
the scale function:
\begin{equation*}
\phi(t)=\sqrt{2}\sum_{k=0}^{M}h(k)\phi(2t-k).
\end{equation*}
There are better methods one could apply for the initialization,
as suggested by Beylkin {\em et al.} in
 \cite{kn:coifman}.

\begin{remark}
We use only the vanishing moments property (\ref{vmoments}) and
compactness of the support of the wavelets to prove Theorems
\ref{mainthm2} and \ref{mainthm1}. However we need the
orthogonality introduced by the multiresolution analysis for
implementing the pyramidal algorithm introduced in Section
\ref{pyramidalalgorithm}. Moreover in our analysis of the S\&P 500
index data we work with compactly supported wavelets to further
increase the efficiency of the pyramidal algorithm.
\end{remark}

\small{
}

\begin{thebibliography}{99}
\bibitem{kn:veitch1} ABRY, P. , and D. VEITCH (1998):
Wavelet Analysis of Long-range-dependent Traffic, \emph{IEEE
Transactions on Information Theory} 44, 2-15.
\bibitem{kn:abry} ABRY, P., and F. SELLAN (1996): The
Wavelet-Based Synthesis for Fractional Brownian Motion Proposed by
F. Sellan and Y. Meyer: Remarks and Fast Implementation,
\emph{Applied and Computational Harmonic Analysis} 3, 377-383
\bibitem{kn:andersen} ANDERSEN, T. G. , and T. BOLLERSLEV (1997):
Heterogeneous Information Arrivals and Return Volatilty Dynamics:
Uncovering the Long-Run in High Frequency Returns, \emph{Journal
of Finance} 52, 975-1005.
\bibitem{kn:andersen2}ANDERSEN, T. G. , and T. BOLLERSLEV (1997):
Intraday Periodicity and Volatility Persistence in Financial
Markets, \emph{Journal of Empirical Finance} 4, 115-158
\bibitem{kn:andersen3}ANDERSEN, T. G. , and T. BOLLERSLEV (1998):
Deutsche Mark-Dollar Volatility: Intraday Activity Patterns,
Macroeconomic Annauncements, and Longer Run Dependencies, \emph{
Journal of Finance} 53, 219-265.
\bibitem{kn:arneodo} ARNEODO, A. (1996): Wavelet Analysis of
Fractals, \emph{Wavelets: Theory and Applications}, Oxford
University Press, New York.
\bibitem{kn:barucci} BARUCCI, E., P. MALLIAVIN, M. E. MANCINO,
R. RENO (2003): The Price-Volatility Feedback Rate: An
Implementable Mathematical Indicator of Market Stability,
\emph{Mathematical Finance}, 13, 17-37.
\bibitem{kn:erhan}  BAYRAKTAR, E. and  H. V. POOR  (2001):
Arbitrage in Fractal Modulated Markets When the Volatility is
Stochastic, \emph{preprint, Princeton University}.
\bibitem{kn:bp} BAYRAKTAR, E., H. V. POOR (2002):
Stochastic Differential Games in a Non-Markovian Setting,
\emph{preprint}, Princeton University.
\bibitem{kn:bhs} BAYRAKTAR, E., U. HORST, K. R. SIRCAR
(2003) A Limit Theorem for Financial Markets with Inert Investors,
\emph{preprint}, Princeton University.
\bibitem{kn:coifman} BEYLKIN, G., R. R. COIFMAN, and V. ROKHLIN
(1992): Wavelets in Numerical Analysis, \emph{Wavelets and Their
Applications}, 181-210. Jones and Barlett.
\bibitem{bjork} BJ\"{O}RK, T. and H. HULT (2003): A Note on the
Self-Financing Condition for the Fractional Black-Scholes Model,
preprint, University of Stockholm.
\bibitem{kn:blackscholes} BLACK, F. and M. SCHOLES (1973): The Pricing of Options and Corporate
Liabilities, \emph{ J. Political Econ.}, 81, 637-659.

\bibitem{kn:cheridito} CHERIDITO, P. (2003): Arbitrage in Fractional Brownian Motion
Models, \emph{Finance and Stochastics}, 7, 533-553.
\bibitem{kn:crr} COX, J., S. ROSS, and M. RUBINSTEIN (1979): Option Pricing: A Simplified
Approach, \emph{J. Financial Economics} 7, 229-263
\bibitem{kn:cutland}CUTLAND, N. J., P. E. KOPP, and W.
WILLINGER (1995): Stock Price Returns and the Joseph Effect: A
Fractal Version of the Black-Scholes model. \emph{Progress in
Probability} 36, 327-351.
\bibitem{kn:icd} DAUBECHIES, I. C. (1992): \emph{Ten Lectures on
Wavelets}, SIAM, Philadelphia.
\bibitem{kn:duncan} {\sc DUNCAN, T. E., Y. HU and B.
PASIK-DUNCAN} (2000): Stochastic Calculus for Fractional Brownian
Motion, \emph{SIAM Journal of Control and Optimization}, 38,
582-612.
\bibitem{elliot} ELLIOTT R. J. and J. VAN DER HOEK (2003): A
General Fractional White Noise Theory and Applications to Finance,
\emph{Mathematical Finance}, 13, 301-330.
\bibitem{kn:fama} FAMA, E. F. (1965): The Behaviour of Stock
Market Prices, \emph{Journal of Business}, 38, 34-105.
\bibitem{kn:fielitz}  GREENE, M. T. and B. D. FIELITZ (1977):
Long-Term Dependence in Common Stock Returns, \emph{Journal of
Financial Economics}, 4, 339-349.
\bibitem{kn:fouque} FOUQUE, J. P., G. PAPANICOLAOU, R. SIRCAR, and K.
SOLNA(2003): Short Time-scale in S\&P 500 Volatility,
\emph{Journal of Computational Finance}, 6(4), 1-23.
\bibitem{kn:rambo2}GEN\c{C}AY, R., F. SEL\c{C}UK, and B. WHITCHER
(2001): Differentiating Intraday Seasonalities Through Wavelet
Multi-Scaling \emph{Physica A} 289, 543-556.
\bibitem{kleinow} HALL, P., W. H\"{A}RDLE, T. KLEINOW and  P.
SCHMIDT (2000): Semiparametric
 Bootstrap Approach To Hypothesis Tests And Confidence Intervals For The Hurst Coefficient,
 \emph{Statistical Inference for Stochastic Processes}, 3,
 263-276.
\bibitem{kn:hida} HIDA, T., H. KUO, J. POTHOFF and L. STREIT (1993):
\emph{White Noise, an Infinite Dimensional Calculus}, Kluwer
Academic Publishers, Boston.
\bibitem{kn:karat} KARATZAS, I. and S.E. SHREVE (1991):
\emph{Brownian Motion and Stochastic Calculus}, Springer-Verlag,
New York.
\bibitem{kn:lo} LO, A. W. (1991): Long-term Memory in Stock Market
Prices, \emph{Econometrica}, 59, 1279-1313.
\bibitem{kn:mallat}MALLAT, S. (2001): \emph{A Wavelet Tour of Signal
Processing}, Academic Press, San Diego.
\bibitem{kn:mandelbrot2} MANDELBROT, B. B. (1963): The
Variation of Certain Speculative Prices, \emph{The Journal of
Business}, 36, 394-419.
\bibitem{kn:mandellong} MANDELBROT, B.B. (1967): Forecasts of
Future Prices, Unbiased Markets and Martingale Models,
\emph{Journal of Business}, 39, 242-255.
\bibitem{kn:mandel} MANDELBROT, B.B. (1971): When Can Price Be
Arbitraged Efficiently? A Limit to the Validity of the Random Walk
and Martingale Models, \emph{Rev. Econom. Statis.}, 53, 225-236.
\bibitem{kn:mandelbrot} MANDELBROT, B. B. (1997): \emph{Fractals
and Scaling in Finance}, Springer-Verlag, New York.
\bibitem{kn:merton69}MERTON, R. C. (1969): Lifetime Portfolio Selection Under Uncertainty: The
Continous-time Case, \emph{Rev. Econom. Statist.} 51, 247-257.

\bibitem{kn:strang} NGUYEN T. and G. STRANG (1997):
\emph{Wavelets and Filter Banks}, Wellesley-Cambridge Press,
Wellesley, MA.
\bibitem{kn:oksendal2} \O KSENDAL, B. and Y. HU (2000):
Fractional White Noise and Applications to Finance, \emph{Infinite
Dimensional Analysis, Quantum Probability and Related Topics}, 6,
1-32.
\bibitem{kn:solna} PAPANICOLAOU, G. and K. SOLNA (2001):
Wavelet Based Estimation of Local Kolmogorov Turbulence, \emph{In
'Long-Range Dependence Theory and Applications', edited by
P.Doukhan, G.Oppenheim, M.S. Taqqu}, Birkhauser, Boston.
\bibitem{kn:peters1} PETERS, E. E. (1991): \emph{Chaos and Order in the Capital
Markets}, John Wiley and Sons, New York.
\bibitem{kn:peters}PETERS, E. E (1994): \emph{Fractal Market Analysis: Applying Chaos Theory
to Investment and Economics}, John Wiley and Sons, New York.
\bibitem{kn:Pro}PROTTER, P. (1990): \emph{Stochastic Integration and Differential
Equations}. Springer-Verlag, Berlin.
\bibitem{kn:rogers}ROGERS, C. (1997): Arbitrage with fractional Brownian
motion, \emph{Mathematical Finance} 7, 95-105.
\bibitem{kn:ruzmaikina} RUZMAIKINA, A. A. (1999):
\emph{Stochastic Calculus with Fractional Brownian Motion}, Ph.D
Dissertation, Department of Physics, Princeton University,
Princeton, NJ.
\bibitem{kn:taqqu} SAMORODNITSKY, G. and M. S. TAQQU (1994):
\emph{Stable Non-Gaussian Random Processes}, Chapman \& Hall,
London.
\bibitem{kn:shiryaev}SHIRYAEV, A. N. (1999): \emph{Essentials of Stochastic
Finance}. World Scientific, Singapore.
\bibitem{valkeila} SOTTINEN, T. and E. VALKEILA (2003): On arbitrage and replication in the fractional Black-Scholes pricing
model, \emph{Statistics and Decisions}, 21, 93-107.
\bibitem{kn:stroock}STROOCK, D. W.(1993): \emph{Probability
Theory}, Cambridge University Press, Cambridge UK.
\bibitem{kn:taqqu2} TEVEROVSKY, V., M. S. TAQQU and W.
WILLINGER (1999): A Critical Look at Lo's Modified R/S Statistic,
\emph{Journal of Statistical Planning and Inference}, 80, 211-227.
\bibitem{kn:taqqu1} WILLINGER, W., M. S. TAQQU and V.
TEVEROVSKY (1999): Stock Market Prices and Long-range Dependence,
\emph{Finance and Stochastics}, 3, 1-13.
\bibitem{kn:zahle} Z\"{A}HLE, M. (1998): Integration with Respect
to Fractal Functions and Stochastic Calculus. I, \emph{Probability
Theory and Related Fields}, 111, 333-374.
\end{thebibliography}
\end{document}